\def\timestamp{%
Time-stamp: <zimna-skola-2002.tex: Friday 03-05-2002 at 13:46:17 (cest)>}
\def\stripname Time-stamp: <#1 #2>{#2}
\edef\filedate{\expandafter\stripname\timestamp}

\documentclass[a4paper]%
              {amsart2000}

\renewcommand\newsymbol[5]{%
\DeclareMathSymbol#1{#3}{\ifcase #2\or AMSa\or AMSb\fi}{"#4#5}}
\newcommand\newbbbletter[2]{%
\DeclareMathSymbol#1{0}{AMSb}{`#2}}
\newbbbletter\C         C    
\newbbbletter\D         D    
\newbbbletter\halfline  H    
\newbbbletter\unitint   I    
\newbbbletter\Miod      M    
\newbbbletter\N         N    
\newbbbletter\psarc     P    
\newbbbletter\Q         Q    
\newbbbletter\R         R    
\newbbbletter\Z         Z    
\usepackage[mathscr]{euscript}
\newcommand\Ba{\mathscr}   
\newcommand\Meas{\Ba{M}}   
\newcommand\smallcard{\mathfrak}
\newcommand\bee{\smallcard{b}}    
\newsymbol\restr   1216         
\newsymbol\le      1236         
\newsymbol\ge      123E         
\newsymbol\symmdif 124D         
\newsymbol\sulk    1061         
\newsymbol\notle       230A
\newsymbol\subsetneq   2328
\newsymbol\nsubseteq 232A
\let\diamond\diamondsuit
\newcommand\famsets{\mathcal}
\newcommand\ideal{\mathcal}
\newcommand\famA{\famsets{A}}
\newcommand\famB{\famsets{B}}
\newcommand\famC{\famsets{C}}
\newcommand\famD{\famsets{D}}
\newcommand\famF{\famsets{F}}
\newcommand\famL{\famsets{L}}
\newcommand\famM{\famsets{M}}
\newcommand\famT{\famsets{T}}
\newcommand\pow {\famsets{P}}  
\newcommand\Null{\ideal{N}}    
\DeclareMathAlphabet{\Bi}{OML}{cmm}{b}{it}
\newcommand{\bH}{\Bi{H}}
\newcommand{\bN}{\Bi{N}}
\let\savephi\phi
\let\phi\varphi
\let\varphi\savephi
\newcommand\cl{\operatorname{cl}}
\newcommand\Int{\operatorname{int}}
\newcommand\Bor{\operatorname{Bor}}

\newcommand\concat{\mathord{{^\sulk}}}
\newcommand\fin{\mathit{fin}}
\newcommand\finsets{[\omega]^{<\omega}}
\newcommand\functions{{}^\omega\omega}
\newcommand\bintree{2^{<\omega}}
\newcommand\sbintree[1][f]{2^{#1(n)}}
\newcommand\almeq{\mathrel{=^*}}
\newcommand\almneq{\mathrel{\neq^*}}
\newcommand\almsubseteq{\mathrel{\subseteq^*}}
\newcommand\alml{\mathrel{<^*}} 
\newcommand\almle{\mathrel{\le^*}}
\newcommand\preim{^\leftarrow} 
\newcommand\powN{\pow(\N)}
\newcommand\powNfin{\powN/\fin}

\newcommand{\Dstar}{\D^*}
\newcommand{\Hstar}{\halfline\mskip1mu{}^*}
\newcommand{\Mstar}{\Miod^*}
\newcommand{\Nstar}{\N^*}
\newcommand\orpr[2]{\langle#1,#2\rangle}
\newcommand\card[1]{\lvert#1\rvert}

\newcommand\axiom[1]{\mathsf{#1}}
\newcommand\CH{\axiom{CH}}
\newcommand\OCA{\axiom{OCA}}
\newcommand\PFA{\axiom{PFA}}

\newcommand\Cantor{2^\omega}
\newcommand\szer{s\concat0}
\newcommand\sone{s\concat1}

\def\iff/{if\kern0ptf}
\newcommand\Hyp[1]{2^{#1}}	
\newcommand\HypX{\Hyp{X}}
\newcommand\Cs[1]{\famC(#1)}	
\newcommand\CsX{\Cs{X}}
\newcommand\Co[1]{\operatorname{Co}{#1}}
\let\meet\wedge
\let\join\vee
\let\bigmeet\bigwedge

\newcommand\0{\mathbf{0}}
\newcommand\1{\mathbf{1}}
\newcommand\cs[1]{\overline{#1}}
\newcommand\conn{\operatorname{conn}}
\let\implies\rightarrow
\DeclareMathSymbol{\lor}{\mathrel}{symbols}{"74}  
\DeclareMathSymbol{\land}{\mathrel}{symbols}{"75} 
\newcommand\LL{\mathfrak{L}}	
\let\elsub\prec
\newtheorem{theorem}{Theorem}[section]
\newtheorem{proposition}[theorem]{Proposition}
\newtheorem{lemma}[theorem]{Lemma}
\theoremstyle{definition}

\theoremstyle{remark}
\newtheorem{remark}[theorem]{Remark}
\theoremstyle{definition}
\newtheorem{definition}[theorem]{Definition}
\newcounter{claim}[theorem]
\renewcommand\theclaim{\arabic{claim}}
\newenvironment{claim}%
           {\endgraf\stepcounter{claim}\smallskip
            \def\proof{\endgraf\vskip1pt plus1pt 
                       \noindent\textsc{Proof}. }
            \noindent\textsc{Claim \theclaim}. \ignorespaces}%
           {\parfillskip0pt\hfil$\symmdif$\endgraf\smallskip}

\begin{document}

\title{An Algebraic and Logical approach to continuous images}

\author{Klaas Pieter Hart}
\address{Faculty of Information Technology and Systems\\
         TU Delft\\
         Postbus 5031\\
         2600~GA {} Delft\\
         the Netherlands}
\email{k.p.hart@its.tudelft.nl}
\urladdr{http://aw.twi.tudelft.nl/\~{}hart}
\date{Zimn\'a \v{S}kola 2002 (\filedate)}
\subjclass[2000]{03C50, 03C98, 03E55, 03E65, 06D05, 06E15, 28A60, 
                 54A35, 54C10, 54D40, 54D80, 54F15, 54F45, 54H10}

\keywords{lattice, ultrafilter, Wallman representation, Wallman space,
          Boolean algebra, Stone space, duality, continuous surjection,
          embedding, universal continuum, \v{C}ech-Stone remainder of~$\R$,
          }

\begin{abstract}
Continuous mappings between compact Hausdorff spaces can be studied using 
homomorphisms between algebraic structures (lattices, Boolean algebras) 
associated with the spaces.
This gives us more tools with which to tackle problems about these continuous
mappings --- also tools from Model Theory.
We illustrate by showing that 1)~the \v{C}ech-Stone remainder $[0,\infty)$ has
a universality property akin to that of~$\Nstar$;
2)~a theorem of Ma\'ckowiak and Tymchatyn implies it own generalization to
non-metric continua; and
3)~certain concrete compact spaces need not be continuous images of~$\Nstar$.
\end{abstract}

\maketitle

\section*{Introduction}

These notes reflect a series of lectures given at the 
\textsl{30th Winterschool on Abstract Analysis} (Section Topology).
In it I surveyed results from the papers~\cite{DowHart2001a}, 
\cite{HartvanMillPol2001}, \cite{DowHart99} and~\cite{DowHart2000a}.
These results are of a topological nature but their proofs involve
algebraic structures associated with the spaces in question.
These proofs also have logical components.
In Sections~\ref{sec.universal} and~\ref{sec.HI}
I use notions from Model Theory show the existence of certain
continua and mappings between them.
In Section~\ref{sec.OCA} we see how the Open Colouring Axiom
implies that very concrete spaces are not continuous images
of~$\Nstar$.

To make these notes reasonably self-contained I devote two sections
to some model-theoretic and algebraic preliminaries.

\section{Lattices and Boolean algebras}
\label{sec.Ba.lattice}

In~\cite{Wallman38} Wallman generalized Stone's representation
theorem for Boolean algebras, from~\cite{StoneMH37a}, to the class
of distributive lattices.
Given a distributive lattice~$L$, with~$\0$ and~$\1$ and operations~$\meet$ 
and~$\join$, we say that $F\subseteq L$ is a \emph{filter} if it
satisfies $\0\notin F$, 
if $a,b\in F$ then $a\meet b\in F$, and
if $a\in F$ and $b\ge a$ then $b\in F$;
as always, an \emph{ultrafilter} is a maximal filter.

The \emph{Wallman representation} (or \emph{Wallman space}) $wL$ 
of~$L$ is the space with the set of all ultrafilters on~$L$ as its
underlying set.
For every~$a\in L$ we put $\cs a=\{u\in wL:a\in u\}$ and we
use the family $\famA=\{\cs a:a\in L\}$ as a base for the closed
sets of a topology on~$wL$.
The resulting space~$wL$ is a compact $T_1$-space and the map
$a\mapsto\cs{a}$ is a homomorphism from~$L$ onto~$\famA$.
The homomorphism is an isomorphism if and only if $L$ is \emph{disjunctive} 
or \emph{separative},
which means: if $a\notle b$ then there is $c\in L$ such that 
$c\le a$ and $c\meet b=\0$.

Every compact $T_1$-space~$X$ can be obtained in this way: $X$~is the Wallman
representation of its own family of closed sets.
From this it is clear that $wL$~is not automatically Hausdorff; in fact
$wL$~is Hausdorff if and only if $L$~is \emph{normal}, which is expressed
as follows:
\begin{equation}\label{L.is.normal}
(\forall x)(\forall y)(\exists u)(\exists v)
\bigl[(x\meet y=\0)\implies
\bigl((x\meet u=\0)\land(y\meet v=\0)\land(u\join v=\1)\bigr)\bigr].
\end{equation}
Note: in order to avoid confusion we write, for the nonce,
$\land$ and $\lor$ for logical `and' and `or' respectively.

The duality is not perfect; one space can represent many different
lattices: one has $X=w\famA$ whenever $\famA$ is a base for the
closed sets of~$X$ that is closed under finite unions and 
intersections --- such a base is referred to as a 
\emph{lattice base for the closed sets}.
For example, the unit interval~$[0,1]$ can also be obtained as the
representation of the lattice consisting of all finite unions
of closed intervals with rational end~points.

Many properties of the space can be expressed using the elements
of~$L$ only.
For example the formula expresses that $wL$~is connected:
\begin{equation}\label{L.is.connected}
(\forall x)(\forall y)
\bigl[\bigl((x\meet y=\0)\land(x\join y=\1)\bigr)\implies
      \bigl((x=\0)\lor(x=\1)\bigr)\bigr].
\end{equation}
This suffices because every lattice base for the closed sets of a compact space
contains every clopen set of that space.
For later use we interpret formula as expressing a property of~$\1$,
to wit ``$\1$ is connected''; we therefore abbreviate it 
as~$\conn(\1)$ and we shall write $\conn(a)$ to denote 
Formula~\ref{L.is.connected} with $\1$ replaced by~$a$ and use it 
to express that $a$~is connected
(or better: the set represented by~$a$ is connected).

\subsection{Boolean algebras}

If $L$~is a Boolean algebra then the family $\{\cs{a}:a\in L\}$
consists exactly of the clopen subsets of~$wL$ and so the space
$wL$ is \emph{zero-dimensional}. 
This makes for a prefect duality between Boolean algebras and 
compact zero-dimensional spaces because every compact
zero-dimensional space represents its own family~$\Co{X}$ of clopen
sets.
This Stone's representation theorem for Boolean algebras;
we call $wL$ the \emph{Stone space} of~$L$.

\subsection{Making continuous surjections}

We use the algebraic approach in the construction of continuous
onto mappings.
The following lemma tells us how this works.

\begin{lemma}\label{lemma.onto.map.lattice}
Let $X$ be compact Hausdorff and $L$ some normal, distributive and disjunctive
lattice.
If $X$ has a lattice base~$\famB$ for the closed sets that is embeddable
into~$L$ then $wL$ admits a continuous surjection onto~$X$.
\end{lemma}

\begin{proof}
We only sketch the argument.
Let $\phi:\famB\to L$ be an embedding and define $f:wL\to X$ by
``$f(p)$~is the unique point in 
$\bigcap\bigl\{C\in\famB:p\in\phi(C)\bigr\}$''.
It is straightforward to check that $f$~is onto and that
$f\preim[C]=\phi(C)$ for all~$C$.
\end{proof}

\section{Elements of Model Theory}
\label{sec.teorie.modelu}

In this section I review some notions and results from Model Theory that we
shall use later on.
Our basic reference for model theory is Hodges' book \cite{Hodges93}.
Dobr\'ym \'uvodem do Teorie Model\accent23u je Kapitola~V
v knize A.~Sochora~\cite{Sochor2001}.

Model Theory studies structures from a logical point of view.
These structures can be groups, fields, ordered sets and,
important for us, lattices.
In what follows I shall try to illustrate the Model-theoretic notions using
ordered sets of fields.

The key notions are those of a \emph{language} and a \emph{theory}.

\subsubsection{Language}

Our languages have two parts.
There is a fixed logical part, consisting of the familiar logical
symbols like $\forall$, $\exists$, $\land$, $\lor$, $\implies$,
$\lnot$, $=$, together with an infinite set of variables.

The second part is specific to the kind of structure that we want
to study.
For example, to study ordered sets we need~$<$;
to study fields we need $+$, $\times$, $0$, and~$1$.

\subsubsection{Theories}

A theory is a set of formulas; nothing more, nothing less.
An interesting theory should be about something non-trivial 
(which may a matter of taste) and consistent, which means
that you cannot derive a false statement from it.

One normally specifies a theory by listing a few formulas as its
starting point (as its axioms) and, tacitly, assumes that the
consequences of these axioms make up the full theory.

Thus, the theory of (linearly) ordered sets has the following three
formulas as its axioms:
\begin{enumerate}
\item $(\forall x)\lnot(x<x)$,
\item $(\forall x)(\forall y)\bigl((x<y)\lor(y<x)\lor(x=y)\bigr)$, and
\item $(\forall x)(\forall y)(\forall z)
      \bigl(((x<y)\land(y<z))\implies (x<z)\bigr)$
\end{enumerate}

Finally then, a \emph{model} for a theory is a structure for the language
where all the formulas of the theory are valid. 
Somewhat tautologically then a linearly ordered set is a model for
the theory of linearly ordered sets.

\subsection{Compactness and completeness}

Two very important theorems, for us, are the Compactness Theorem
and the Completeness Theorem.

The Compactness Theorem says that a theory is consistent if and only
if every finite subset is consistent.
Given the definition of consistency that we adopted this is a 
actually a triviality: any derivation uses only a finite set of 
formulas.
The Compactness Theorem gets quite powerful in combination with the
Completeness Theorem, which says that a theory is consistent if and
only if it has a model.
We shall use the nontrivial consequence that a theory has a model
if and only if every finite subset has a model.

\subsection{Elementarity}

Two structures are \emph{elementarily equivalent} if they satisfy the same
sentences (formulas without free variables); 
informally: they look superficially the same.
For example the ordered sets $\Q$, $\N$ and $\Z$ are all different:
consider the sentences 
$(\forall x)(\forall y)(\exists z)\bigl((x<y)\implies(x<z)\land(z<y)\bigr)$;
$(\exists x)(\forall y)(x\le y)$ and
$(\exists x)(\exists y)(\forall z)\bigl((z\le x)\lor(y\le z)\bigr)$.

On the other hand the ordered sets $\Q$ and $\R$ are elementarily
equivalent.
This can be gleaned from the material in Chapter~3 of~\cite{Hodges93}.
For us it is important to know that any two atomless Boolean 
algebras are elementarily equivalent
\cite[Theorem~5.5.10]{ChangKeisler77}.

\subsubsection{Elementary substructures}

We say that $A$ is an \emph{elementary substructure} of $B$,
written $A\elsub B$, if every equation with parameters in~$A$
that has a solution in~$B$ already has a solution in~$A$.

The field~$\Q$ is not an elementary substructure of the 
field~$\R$, consider the equation $x^2=2$.
On the other hand, the field of algebraic numbers is an elementary
substructure of the field~$\C$ of complex numbers
(see Appendix~A.5 of~\cite{Hodges93}).

The L\"owenheim-Skolem theorem provides us with many elementary
substructures:
if $A$ is a structure for a language~$\LL$ and $X\subseteq A$ then
there is an elementary substructure~$B$ of~$A$ with $X\subseteq B$
and $\card{B}\le\card{X}\cdot\card{\LL}\cdot\aleph_0$.
Normally the language $\LL$ is countable, so that we can get
many \emph{countable} elementary substructures; we will use 
this often to construct metric continua.

\subsection{Saturation}

Given a cardinal $\kappa$ one calls a structure (e.g., a field, 
a group, an ordered set, a lattice) is said to be 
$\kappa$-saturated if, loosely speaking, every consistent set of 
equations, of cardinality less than~$\kappa$ and with parameters
from the given structure, has a solution, 
where a set of equations is consistent if every finite subsystem 
has a solution 
\emph{possibly at first in some extension of the given structure}.
Thus, e.g., $\{{0<z},{z<1}\}$ is consistent in~$\N$, because
a solution can be found in the extension $\N\cup\{\frac12\}$; 
on the other hand $\{{z<0},{1<z}\}$ is clearly inconsistent.
As the first system has no solution in~$\N$ itself it witnesses
that $\N$~is not $\aleph_0$-saturated.

Going one step up, the ordered set of the reals is \emph{not} 
$\aleph_1$-saturated because the following countable system of
equations, though consistent, does not have a solution: 
$0<x$ together with $x<\frac1n$ ($n\in\N$).
On the other hand, any ultrapower~$\R^\omega_u$ of~$\R$ is
$\aleph_1$-saturated as an ordered set --- see \cite[Theorem~9.5.4]{Hodges93}.
Such an ultrapower is obtained by taking the power~$\R^\omega$,
an ultrafilter~$u$ on~$\omega$ and identifying points~$x$ and~$y$
if $\{n:x_n=y_n\}$ belongs to~$u$.
The ordering~$<$ is defined in the obvious way:
$x<y$ iff $\{n:x_n<y_n\}$ belongs to~$u$.
It is relatively easy to show that this gives an $\aleph_1$-saturated
ordering; given a countable consistent set of equations $x<a_i$ and $x>b_i$
($i\in\omega$), one has to produce a single~$x$ that satisfies them all;
the desired~$x$ can be constructed by a straightforward diagonalization.

\subsection{Universality}

Finally, a structure is $\kappa$-universal if it contains a copy of every
structure of cardinality less than~$\kappa$ that is elementarily equivalent 
to it.

Our last ingredient is Theorem~10.1.6 from~\cite{Hodges93}, which
states that $\kappa$-saturated structures are $\kappa^+$-universal.
When we apply this to an ultrapower~$\R^\omega_u$ then we find that
it contains an isomorphic copy of every $\aleph_1$-sized dense 
linear order without end~points --- a result that can also be
established directly by a straightforward transfinite recursion.
It also follows that $\R^\omega_u$ contains an isomorphic copy of
\emph{every} $\aleph_1$-sized linear order: simply make it dense by inserting a
copy of the rationals between any pair of neighbours and attach copies of the
rationals at the beginning and the end to get rid of possible end~points; the
resulting ordered set is still of cardinality~$\aleph_1$ and can therefore be
embedded into~$\R^\omega_u$.

\section{Universal compact spaces}
\label{sec.universal}

Here we combine the algebra and model theory to provide proofs of
universality of certain spaces.
Here `universality' is meant in the mapping-onto sense, i.e.,
space~$X$ is universal for a class of spaces if it belongs to
the class and every space in the class is a continuous image of~$X$.

\subsection{The Cantor set and $\Nstar$}

Let us begin by reviewing two well-known theorems from topology.
The first is due to Alexandroff~\cite{Alexandroff27} and 
Hausdorff~\cite{Hausdorff27}; it states that every
compact metric space is a continuous image of the Cantor set.
The second is Parovi\v{c}enko's theorem \cite{Parovicenko63}
that every compact space of weight~$\aleph_1$ (or less) is a continuous image 
of the space~$\Nstar$.
Both theorems can be proven in a similar fashion.
The first step is a theorem of Alexandroff~\cite{Alexandroff36}.

\begin{theorem}
Every compact Hausdorff space is the continuous image of a compact
zero-dimensional space of the same weight.
\end{theorem}

\begin{proof}
Let $\famB$ be a base for the space~$X$, of size~$w(X)$.
Let $\Ba{B}$ be the Boolean subalgebra of~$\pow(X)$ generated
by~$\famB$.
The Stone space~$Y$ of~$\Ba{B}$ is the sought-after space.
If $u\in Y$ (so $u$~is an ultrafilter on~$\Ba{B}$) then
$\bigcap\{\cl B:B\in u\}$ consists of one point~$x_u$; the map
$u\mapsto x_u$ is a continuous from~$Y$ onto~$X$.
\end{proof}

The second step is to embed the clopen algebra of~$Y$, which happens
to be~$\Ba{B}$, into the clopen algebra of the Cantor set 
or~$\Nstar$ respectively --- the Lemma~\ref{lemma.onto.map.lattice}
applies to give a continuous map from the Cantor set (or~$\Nstar$)
onto~$Y$.
We do this in a roundabout way, to set the stage for a similar proof
involving continua.
First we embed $\Ba{B}$ into the clopen algebra~$\Ba{C}$ 
of~$Y\times\Cantor$ (in the obvious way), this latter algebra is
atomless.

It is fairly straightforward to show that atomless Boolean algebras
are $\aleph_0$-sat\-u\-rated and it is a little more work to show that
the clopen algebra of~$\Nstar$ (which is $\powNfin$) is 
$\aleph_1$-saturated (see~\cite{JonssonOlin68}).

We see that every countable atomless Boolean algebra is embeddable
into the clopen algebra of~$\Cantor$ and every atomless Boolean
algebra of size~$\aleph_1$ (or less) is embeddable into~$\powNfin$.
But this exactly what we still needed to establish.

\subsection{A universal continuum}

In this section we shall apply the ideas developed above in a proof
that the \v{C}ech-Stone remainder of~$[0,\infty)$ maps onto every 
continuum of weight~$\aleph_1$ or less.

\subsubsection{The continuum $\Hstar$}

We write $\halfline=[0,\infty)$ and show that the continuum~$\Hstar$
maps onto every continuum of weight~$\aleph_1$.
This continuum has a nice base for its closed sets: the lattice
$$
\famL=\{A^*:\text{$A$~is closed in~$\halfline$}\}.
$$
Here, as is common, $A^*$ abbreviates $\cl A\cap\Hstar$.
Another way to represent this lattice is as the quotient of the
lattice~$\Hyp\halfline$ by the ideal of compact sets.
Therefore one way to apply Lemma~\ref{lemma.onto.map.lattice}
would be to construct, given a continuum~$X$ of weight~$\aleph_1$
or less,
a lattice base~$\famB$ for the closed sets of~$X$ and a map
$\phi:\famB\to\Hyp\halfline$ whose composition with the quotient
homomorphism is a lattice embedding.
Unfortunately this does not seem to be easy to do, even for metric
continua.

\subsubsection{The metric case}

Our starting point is the following theorem,
due to Aarts and van Emde Boas \cite{AartsvanEmdeBoas67};
as we shall need this theorem and its proof later, 
we provide a short argument.

\begin{theorem}\label{thm.jan}
The space\/ $\Hstar$ maps onto every metric continuum.
\end{theorem}

\begin{proof}
Consider a metric continuum~$K$ and assume it is embedded into
the Hilbert cube~$Q=[0,1]^\infty$.
Choose a countable dense subset~$A$ of\/~$K$ and enumerate it
as~$\{a_n:n\in\omega\}$.
Next choose, for every~$n$, a finite sequence of points
~$a_n=a_{n,0}$, $a_{n,1}$, \dots,~$a_{n,k_n}=a_{n+1}$
such that $d(a_{n,i},a_{n,i+1})<2^{-n}$ for all~$i$ --- this uses the
connectivity of\/~$K$.
Finally, let~$e$ be the map from~$\halfline$
to~$(0,1]\times Q$ with first coordinate~$e_1(t)=2^{-t}$
and whose second coordinate satisfies $e_2(n+\frac i{k_n})=a_{n,i}$
for all~$n$ and~$i$ and is (piecewise) linear otherwise.

It is clear that $e$ is an embedding, and one readily checks that
$\cl e[\halfline]=e[\halfline]\cup\bigl(\{0\}\times K\bigr)$;
the \v{C}ech-Stone extension~$\beta e$ of\/~$e$ maps~$\Hstar$ onto~$K$.
\end{proof}

This theorem and its proof give us an almost lattice-embedding for
bases of metric continua.

\begin{lemma}\label{lemma.first.step}
Let $K$ be a metric continuum and let $x\in K$.
There is a map~$\phi$ from~$\Hyp{K}$ to~$\Hyp{\halfline}$ such that
\begin{enumerate}
\item $\phi(\emptyset)=\emptyset$ and $\phi(K)=\halfline$;\label{fs.i}
\item $\phi(F\cup G)=\phi(F)\cup\phi(G)$;\label{fs.ii}
\item if $F_1\cap\cdots\cap F_n=\emptyset$\label{fs.iii}
      then $\phi(F_1)\cap\cdots\cap\phi(F_n)$ is compact; and
\item $\N\subseteq\phi\bigl(\{x\}\bigr)$.\label{fs.iv}
\end{enumerate}
In addition, if some countable family~$\famC$ of nonempty closed subsets
of\/~$K$ is given in advance, then we can arrange that for every~$F$ in~$\famC$
the set~$\phi(F)$ is not compact.
\end{lemma}

\begin{proof}
As proved in Theorem~\ref{thm.jan}, there is a map from~$\Hstar$ onto~$K$.

The proof given in~\cite{AartsvanEmdeBoas67} 
(and the one given above) is flexible enough to allow us
to ensure that the embedding~$e$ of\/~$\halfline$ into~$(0,1]\times Q$ is such
that $e(n)=\langle2^{-n},x\rangle$ for every~$n\in\N$
and that for every element~$y$ of some countable set~$C$ the set
$\{t:e_1(t)=y\}$ is cofinal in~$\halfline$ --- it is also easy to change the
description of~$e$ in the proof we gave to produce another~$e$ with the
desired properties.
In our case we let $C$ be a countable subset of\/~$K$ that meets every
element of the family~$\famC$.

We now identify $K$ and $\{0\}\times K$, and define a map
$\psi:\Hyp{K}\to\Hyp{\unitint\times Q}$ by
\[
\psi(F)=\bigl\{y\in:\unitint\times Q:d(y,F)\le d(y,K\setminus F)\bigr\}.
\]
In \cite[\S\,21\,XI]{Kuratowski66} it is shown that for all $F$ and $G$
we have
\begin{itemize}
\item $\psi(F\cup G)=\psi(F)\cup\psi(G)$;
\item $\psi(K)=\unitint\times Q$ and $\psi(\emptyset)=\emptyset$ ---
      by the fact that $d(y,\emptyset)=\infty$ for all~$y$; and
\item $\psi(F)\cap K=F$.
\end{itemize}
Note that for every~$y\in K$ we have
$d\bigl(\langle t,y\rangle,\{y\}\bigr)
 =d\bigl(\langle t,y\rangle,K\setminus\{y\}\bigr)=t$, and hence
$\unitint\times\{y\}\subseteq\psi\bigl(\{y\}\bigr)$.

Now define $\phi(F)=e\preim\bigl[\psi(F)\bigr]$---or rather, after identifying
$\halfline$ and~$e[\halfline]$, set $\phi(F)
\allowbreak
=\psi(F)\cap e[\halfline]$.
All desired properties are easily verified: \ref{fs.i} and \ref{fs.ii}
are immediate;
to see that \ref{fs.iii}~holds, note that
if $F_1\cap\cdots\cap F_n=\emptyset$
then $\cl\phi(F_1)\cap\cdots\cap\cl\phi(F_n)\cap K=\emptyset$,
so that $\cl\phi(F_1)\cap\cdots\cap\cl\phi(F_n)$ is a compact subset
of\/~$\halfline$.
That \ref{fs.iv}~holds follows from the way we chose the values~$e(n)$ for
$n\in\N$.

Finally, if $F\in\famC$ and $y\in C\cap F$, then
the cofinal set $\{t:\pi\bigl(e(t)\bigr)=y\}$ is a subset of\/~$\phi(F)$,
so that $\phi(F)$~is not compact.
\end{proof}

\subsubsection{Making continuous surjections (bis)}

Lemma~\ref{lemma.first.step} indicates that 
Lemma~\ref{lemma.onto.map.lattice} may not be directly applicable.
On the other hand, it does indicate that lattice-embeddings may
not be necessary for obtaining onto mappings.
The following theorem shows how much we actually need.

\begin{theorem}\label{thm.how.to.map.onto}
Let $X$ and $Y$ be compact Hausdorff spaces and let $\famC$~be a base for the
closed subsets of\/~$Y$ that is closed under finite unions and finite
intersections.
Then $Y$~is a continuous image of\/~$X$ if and only if there is a map
$\phi:\famC\to\Hyp{X}$ such that
\begin{enumerate}
\item $\phi(\emptyset)=\emptyset$, and\label{cond.i}
      if $F\neq\emptyset$ then $\phi(F)\neq\emptyset$;
\item if $F\cup G=Y$ then $\phi(F)\cup\phi(G)=X$; and\label{cond.ii}
\item if $F_1\cap\cdots\cap F_n=\emptyset$\label{cond.iii}
      then $\phi(F_1)\cap\cdots\cap\phi(F_n)=\emptyset$.
\end{enumerate}
\end{theorem}

\begin{proof}
Necessity is easy: given a continuous onto map $f:X\to Y$, let
$\phi(F)=f\preim[F]$.
Note that $\phi$ is in fact a lattice-embedding.

To prove sufficiency, let $\phi:\famC\to\Hyp{X}$ be given and consider for
each $x\in X$ the family $\famF_x=\bigl\{F\in\famC:x\in\phi(F)\bigr\}$. We
claim that $\bigcap\famF_x$ consists of exactly one point. Indeed, by
condition~\ref{cond.iii} the family~$\famF_x$ has the finite intersection
property, so that $\bigcap\famF_x$ is nonempty. Next assume that $y_1\neq y_2$
in~$Y$ and take $F,G\in\famC$ such that $F\cup G=Y$,\ $y_1\notin F$ and
$y_2\notin G$. Then, by condition~\ref{cond.ii}, either $x\in\phi(F)$ and so
$y_1\notin\bigcap\famF_x$, or $x\in\phi(G)$ and so
$y_2\notin\bigcap\famF_x$.

We define $f(x)$ to be the unique point in~$\bigcap\famF_x$.

To demonstrate that $f$ is continuous and onto, we show that for every closed
subset~$F$ of\/~$Y$ we have
\begin{equation} \tag{$*$}
f\preim[F]=\bigcap\bigl\{\phi(G):G\in\famC, F\subseteq \Int G\bigr\}.
\end{equation}
This will show that preimages of closed sets are closed and that every fiber
$f\preim(y)$~is nonempty.

We first check that the family on the right-hand side has the finite
intersection property.
Even though $F$ and the complement~$K$ of~$\bigcap_i\Int G_i$ need not
belong to~$\famC$, we can still find~$G$ and~$H$ in~$\famC$ such that
$G\cap K=H\cap F=\emptyset$ and $H\cup G=Y$.
Indeed, apply compactness and the fact that $\famC$~is a lattice to
find~$C$ in~$\famC$ such that $F\subseteq C\subseteq\bigcap_i\Int G_i$ and
then $D\in\famC$ with $K\subseteq D$ and $D\cap C=\emptyset$;
then apply normality of~$\famC$ to~$C$ and~$D$.
Once we have $G$ and~$H$ we see that for each~$i$ we also have~$H\cup G_i=Y$,
and so $\phi(H)\cup\phi(G_i)=X$; combined with
$\phi(G)\cap\phi(H)=\emptyset$, this gives
$\phi(G)\subseteq\bigcap_i\phi(G_i)$.

To verify~$(*)$, first let $x\in X\setminus f\preim[F]$.
As above we find $G$ and~$H$ in~$\famC$ such that
$f(x)\notin G$, $H\cup G=Y$ and $H\cap F=\emptyset$.
The first property gives us~$x\notin\phi(G)$; the other two imply
that $F\subseteq\Int G$.

Second, if $F\subseteq\Int G$, then we can find $H\in\famC$ such that $H\cup
G=X$ and $H\cap F=\emptyset$. It follows that if $x\notin\phi(G)$ we have
$x\in\phi(H)$;  hence $f(x)\in H$ and so $f(x)\notin F$.
\end{proof}

We shall now show how to construct, given a continuum $K$ of weight~$\aleph_1$, 
a map~$\phi$ from a base for the closed sets of~$K$ into the 
base~$\famL$ as in Theorem~\ref{thm.how.to.map.onto}.
Our plan is to find this map using the model-theoretic machinery 
described above.

This would require two steps.
Step~1 would be to show that $\famL$~is an $\aleph_1$-saturated
lattice and hence $\aleph_2$-universal.
Step~2 would then be to show that every lattice of size~$\aleph_1$
is embeddable into a lattice of size~$\aleph_1$ that itself is 
elementarily equivalent to~$\famL$.

There are two problems with this approach:
1)~we were not able to show that $\famL$~is
$\aleph_1$-saturated, and
2)~Lemma~\ref{lemma.first.step} does not give a lattice-embedding.
We shall deal with these problems in turn.

\subsubsection{An $\aleph_1$-saturated structure}

As mentioned above, we do not know whether $\famL$~is
$\aleph_1$-saturated.
We can however find an $\aleph_1$-saturated sublattice:
\[
\famL'=\{A^*:\text{$A$~is closed in~$\halfline$, and
           $\N\subseteq A$ or $\N\cap A=\emptyset$}\,\}.
\]
This lattice is a base for the closed sets of the space~$\bH$, obtained
from~$\Hstar$ by identifying~$\Nstar$ to a point.

To see that $\famL'$ is $\aleph_1$-saturated we introduce another space,
namely $\Miod=\omega\times\unitint$, where $\unitint$~denotes the unit
interval.
The canonical base~$\famM$ for the closed sets of~$\Mstar$ is naturally
isomorphic to the reduced power $(2^\unitint)^\omega$ modulo the
cofinite filter.
It is well-known that this structure is $\aleph_1$-saturated
--- see~\cite{JonssonOlin68}.
The following substructure~$\famM'$ is $\aleph_1$-saturated as well:
\[
\famM'=\{A^*:\text{$A$~is closed in~$\Miod$, and
             $\bN\subseteq A$ or
             $\bN\cap A=\emptyset$}\,\},
\]
where $\bN=\{0,1\}\times\omega$.
Indeed, consider a countable consistent set $T$ of equations with constants
from~$\famM'$.
We can then add either $\bN\subseteq x$ or $\bN\cap x=\emptyset$
to~$T$ without losing consistency.
Any element of~$\famM$ that satisfies the expanded~$T$ will automatically
belong to~$\famM'$.

We claim that $\famL'$ and $\famM'$ are isomorphic.
To see this, consider the map $q:\Miod\to\halfline$ defined by
$q(n,x)=n+x$.
The \v{C}ech-Stone extension of~$q$ maps $\Miod^*$ onto~$\Hstar$, and it
is readily verified that $L\mapsto q^{-1}[L]$ is an isomorphism
between~$\famL'$ and~$\famM'$.
(In topological language: the space~$\bH$ is also obtained from~$\Mstar$
by identifying~$\bN^*$ to a point.)

\subsubsection{A new language}

The last point that we have to address is that Lemma~\ref{lemma.first.step}
does not provide a lattice embedding, but rather a map that only partially
preserves unions and intersections.
This is where Theorem~\ref{thm.how.to.map.onto} comes in:
we do not need a full lattice embedding, but only a map that preserves certain
identities.
We abbreviate these identities as follows:
\begin{align*}
J(x,y) &\quad\equiv\quad x\vee y=1 ,\\
M_n(x_1,\ldots,x_n) &\quad\equiv\quad x_1\wedge\cdots\wedge x_n=0.
\end{align*}
We can restate the conclusion in Theorem~\ref{thm.how.to.map.onto} in the
following manner: $Y$~is a continuous image of~$X$ if and only if there is
an~$\LL$-homomorphism from~$\famC$ to~$\Hyp{X}$, where $\LL$~is the language
that has $J$ and the~$M_n$ as its predicates \emph{and} where $J$ and the~$M_n$
are interpreted as above.

Note that by considering a lattice with $0$ and~$1$ as an $\LL$-structure
we do not have to mention $0$ and~$1$ anymore; they are implicit in the
predicates.
For example, we could define a normal $\LL$-structure to be one in which
$M_2(a,b)$ implies
$(\exists c,d)\bigl(M_2(a,d)\land M_2(c,b)\land J(c,d)\bigr)$.
Then a lattice is normal iff it is normal as an $\LL$-structure.

\subsubsection{The proof}

Let $\famC$ be a base of size~$\aleph_1$ for the closed sets of the
continuum~$K$.
We want to find an $\LL$-structure~$\famD$ of size~$\aleph_1$
that contains~$\famC$ and that is elementarily equivalent to~$\famL'$.
To this end we consider the diagram of\/~$\famC$; that is, we add the elements
of\/~$\famC$ to our language~$\LL$ and we consider the set~$D_\famC$ of all
atomic sentences from this expanded language that hold in~$\famC$.
For example, if $a\cap b=\emptyset$ and $c\cup d=K$, then
$M_2(a,b)\land J(c,d)$ belongs to~$D_\famC$.

To $D_\famC$ we add the theory~$T_{\famL'}$ of\/~$\famL'$,
to get a theory~$T_\famC$.
Let $\famC'$ be any countable subset of~$\famC$ and assume, without loss of
generality, that $\famC'$~is a normal sublattice of~$\famC$.
The Wallman space~$X$ of~$\famC'$ is metrizable, because $\famC'$~is countable,
and connected because it is a continuous image of~$K$.
We may now apply Lemma~\ref{lemma.first.step} to obtain an $\LL$-embedding
of~$\famC'$ into~$\famL'$; indeed, condition~4 says that $\Nstar$~will be mapped
onto a fixed point~$x$ of~$X$.

This shows that, for every countable subset~$\famC'$ of~$\famC$, the union
of $D_{\famC'}$ and~$T_{\famL'}$ is consistent, and so, by the compactness
theorem, the theory~$T_\famC$ is consistent.
Let $\famD$ be a model for~$T_\famC$ of size~$\aleph_1$.
This model is as required: it satisfies the same sentences
as~$\famL'$ and it contains a copy of\/~$\famC$, to wit the set of
interpretations of the constants from~$\famC$.

\section{Hereditarily indecomposable continua}
\label{sec.HI}

The model-theoretic approach is also quite useful in the theory of
hereditarily indecomposable continua.

A continuum is \emph{decomposable} if it can be written as the union of two 
proper subcontinua; it is \emph{indecomposable} otherwise.
A \emph{hereditarily indecomposable} continuum is one in which every 
subcontinuum is indecomposable.
It is easily seen that this is equivalent to saying that whenever two
continua in the space meet one is contained in the other.

This latter statement makes sense for arbitrary compact Hausdorff spaces,
connected or not; we therefore extend this definition and call a compact 
Hausdorff space \emph{hereditarily indecomposable} if it satisfies the 
statement above: whenever two continua in the space meet one is contained in 
the other.
Thus, zero-dimensional spaces are hereditarily indecomposable too.

We shall mainly use a characterization of hereditary indecomposability
that can be gleaned from \cite[Theorem~3]{KrasinkiewiczMinc1977} and which
was made explicit in~\cite[Theorem~2]{OversteegenTymchatyn86}.
To formulate it we introduce some terminology.

Let $X$ be compact Hausdorff and let $C$ and $D$ be disjoint closed subsets
of~$X$; as in~\cite{KrasinkiewiczMinc1977} we say that $(X,C,D)$~is
\emph{crooked} between the neighbourhoods $U$ of~$C$ and $V$ of~$D$ if we can 
write $X=X_0\cup X_1\cup X_2$, where each~$X_i$ is closed and, moreover,
$C\subseteq X_0$, $X_0\cap X_1\subseteq V$, $X_0\cap X_2=\emptyset$,
$X_1\cap X_2\subseteq U$ and $D\subseteq X_2$.
We say $X$ is crooked between $C$ and $D$ if $(X,C,D)$~is crooked between
any pair of neighbourhoods of~$C$ and~$D$.

We can now state the characterization of hereditary indecomposability
that we will use.

\begin{theorem}[Krasinkiewicz and Minc]
A compact Hausdorff space is hereditarily indecomposable if and only if it is 
crooked between every pair of disjoint closed (nonempty) subsets.
\end{theorem}

This characterization can be translated into terms of closed sets only;
we simply put $F=X\setminus V$ and $G=X\setminus U$, and reformulate
some of the premises and the conclusions.
We get the following formulation.

\begin{theorem}\label{thm.zig-zag}
A compact Hausdorff space $X$ is hereditarily indecomposable if and only if 
whenever four closed sets $C$, $D$, $F$ and $G$ in~$X$ are given such that 
$C\cap D=C\cap F=G\cap D=\emptyset$ one can write $X$ as the union of three 
closed sets $X_0$, $X_1$ and~$X_2$ such that
$C\subseteq X_0$,
\ $D\subseteq X_2$,
\ $X_0\cap X_1\cap G=\emptyset$,
\ $X_0\cap X_2=\emptyset$, and
  $X_1\cap X_2\cap F=\emptyset$.
\end{theorem}

To avoid having to write down many formulas we call a quadruple
$(C,D,F,G)$ with $C\cap D=C\cap F=D\cap G=\emptyset$ a
\emph{pliand foursome} and we call a triple $(X_0,X_1,X_2)$
with $C\subseteq X_0$,
\ $D\subseteq X_2$,
\ $X_0\cap X_1\cap G=\emptyset$,
\ $X_0\cap X_2=\emptyset$, and
  $X_1\cap X_2\cap F=\emptyset$
a \emph{chicane} for $(C,D,F,G)$.
Thus, a compact Hausdorff space is hereditarily indecomposable if and only if
there is a chicane for every pliand foursome.

This characterization can be improved by taking a base $\famB$ for the closed
sets of the space~$X$ that is closed under finite intersections.
The space is hereditarily indecomposable if and only if there is a chicane
for every pliand foursome whose terms come from~$\famB$.

To prove the nontrivial implication let $(C,D,F,G)$ be a pliand foursome
and let $(O_C,O_D,O_F,O_G)$ be a swelling of it, i.e.,
every $O_P$ is an open set around~$P$ and $O_P\cap O_Q=\emptyset$ if and only
if $P\cap Q=\emptyset$, where $P$ and $Q$ run through $C$, $D$, $F$ and~$G$
(see \cite[7.1.4]{Engelking89}).
Now compactness and the fact that $\famB$ is closed under finite intersections
guarantee that there are $C'$, $D'$, $F'$ and~$G'$ in~$\famB$ such that
$P\subseteq P'\subseteq O_P$ for $P=C$, $D$, $F$,~$G$.
Any chicane for $(C',D',F',G')$ is a chicane for~$(C,D,F,G)$.

\subsection{Hereditarily indecomposable continua of arbitrary weight}

Model theory can help to show that there are hereditarily indecomposable 
continua of arbitrary large weight.
We obtain such continua as Wallman spaces of suitable lattices.

To ensure that $wL$~is hereditarily indecomposable it suffices to have
a chicane for every pliand foursome from~$L$ and this is exactly
what the following formula expresses.
\begin{multline}\label{L.is.HI}
(\forall x)(\forall y)(\forall u)(\forall v)(\exists z_1,z_2,z_3)
\bigl[\bigl((x\meet y=\0)\land(x\meet u=\0)\land(y\meet v=\0)\bigr)\implies{}\\
{}\implies\bigl((x\meet(z_2\join z_3)=\0)\land (y\meet(z_1\join z_2)=\0)
   \land (z_1\meet z_3=\0){}\\
{}\land(z_1\meet z_2\meet v=\0)\land(z_2\meet z_3\meet u=\0)\land
 (z_1\join z_2\join z_3=\1)      \bigr)\bigr].
\end{multline}

The existence of the pseudoarc~$\psarc$ implies that there are one-dimensional
hereditarily indecomposable continua of arbitrarily large weight.
Indeed, the family of closed sets of~$\psarc$ is a distributive and disjunctive
lattice that satisfies formulas~\ref{L.is.normal}, \ref{L.is.connected}
and~\ref{L.is.HI}; it also satisfies
\begin{multline}\label{dimL<=1}
(\forall x_0)(\forall y_0)(\forall x_1)(\forall y_1)
(\exists u_0,v_0,u_1,v_1)
\bigl[\bigl((x_0\meet y_0=\0)\land(x_1\meet y_1=\0)\implies{}\\
{}\implies\bigl((x_0\meet u_0=\0)\land (y_0\meet v_0=\0)\land
                (x_1\meet u_1=\0)\land (y_1\meet v_1=\0)\land{}\\
        {}\land (u_0\join v_0=\1)\land (u_1\join v_1=\1)\land
                (u_0\meet v_0\meet u_1\meet v_1=\0)\bigr)\bigr].
\end{multline}
This formula expresses $\dim wL\le1$ in terms of closed sets; it is
the Theorem on Partitions, see~\cite[Theorem~7.2.15]{Engelking89}.
Therefore this combination of formulas is consistent and so, by the (upward)
L\"owenheim-Skolem theorem, it has models of every cardinality.
Thus, given a cardinal~$\kappa$ there is a distributive and disjunctive
lattice~$L$ of cardinality~$\kappa$ that satisfies~\ref{L.is.normal}, 
\ref{L.is.connected}, \ref{L.is.HI} and~\ref{dimL<=1}.
The space~$wL$ is compact Hausdorff, connected, hereditarily indecomposable,
one-dimensional and of weight~$\kappa$ or less, but with at least~$\kappa$
closed sets.
Thus, if $\kappa\ge2^\lambda$ then the weight of~$wL$ is at least~$\lambda$.

To get a space of weight exactly~$\kappa$ we make sure that $wL$ has at least
$2^\kappa$ many closed sets.
To this end we introduce two sets of $\kappa$~many constants
$\{a_\alpha:\alpha<\kappa\}$ and $\{b_\alpha:\alpha<\kappa\}$ and two sets
of $\kappa$~many formulas:
for every~$\alpha$ the formula $a_\alpha\meet b_\alpha=\0$ and for
any pair of disjoint finite subsets~$p$ and~$q$ of~$\kappa$ the formula
$\bigmeet_{\alpha\in p}a_\alpha\meet\bigmeet_{\alpha\in q}b_\alpha\neq\0$.
Thus we have expanded the language of lattices by a number of constants
and we have added a set of formulas to the formulas that we used above.
This larger set~$\famT_\kappa$ of formulas is still consistent.

Take a finite subset~$T$ of $\famT_\kappa$ and fix a finite subset~$t$ 
of~$\kappa$ such that whenever $a_\alpha\meet b_\alpha=\0$ or
$\bigmeet_{\alpha\in p}a_\alpha\meet\bigmeet_{\alpha\in q}b_\alpha\neq\0$
belong to~$T$ we have $\alpha\in t$ and $p\cup q\subseteq t$.
Now take a map~$f$ from~$\psarc$ onto the cube~$\unitint^t$ and interpret
$a_\alpha$ by~$f\preim[A_\alpha]$ and $b_\alpha$ by~$f\preim[B_\alpha]$;
in this way we have ensured that every formula from~$T$ holds in the family
of closed subsets of~$\psarc$.
Therefore $T$~is a consistent set of formulas and so, because it was arbitrary
and by the compactness theorem, the full set~$\famT_\kappa$ is 
consistent.

Because $\famT_\kappa$ has cardinality~$\kappa$ it has a model~$L$
of cardinality~$\kappa$.
Now $wL$~is as required: its weight is at most~$\kappa$ because $L$~is a base
of cardinality~$\kappa$.
On the other hand: for every subset~$S$ of~$\kappa$, we have, by compactness,
a nonempty closed set
$$
F_S=\bigcap_{\alpha\in S}a_\alpha\cap\bigcap_{\alpha\notin S}b_\alpha
$$
such that $F_S\cap F_T=\emptyset$ whenever $S\neq T$.

\begin{remark}
The reader may enjoy modifying the above argument so as to ensure
that $\bigl\{(a_\alpha,b_\alpha):\alpha<\kappa\bigr\}$ is an essential family
in~$wL$.
To this end write down, for every finite subset~$a$ of~$\kappa$,
a formula~$\phi_a$ that expresses that 
$\bigl\{(a_\alpha,b_\alpha):\alpha\in a\bigr\}$
is essential.
Theorem~2.1 from~\cite{HartvanMillPol2001} more than ensures that the
set of formulas consisting of~\ref{L.is.normal}, \ref{L.is.connected}, 
\ref{L.is.HI} and the~$\phi_a$ is consistent. 
\end{remark}

\subsection{Hereditarily indecomposable preimages}

In~\cite[(19.3)]{MackowiakTymchatyn1984} it is proven that every
metric continuum is the weakly confluent image of some hereditarily 
indecomposable metric curve.
A map is \emph{weakly confluent} if every continuum in the range is the
image of a continuum in the domain.
Using our model-theoretic approach we can generalize this result to
uncountable weights.

\subsubsection{Making an onto map}

To get a (one-dimensional) hereditarily indecomposable
continuum that maps onto the given continuum~$X$ we need to construct a 
distributive, disjunctive and normal lattice~$L$ that satisfies 
formulas~\ref{L.is.connected} and~\ref{L.is.HI} (and~\ref{dimL<=1}), 
and an embedding~$\phi$ of some base~$\famB$ for the closed sets of~$X$ 
into~$L$.

Let a continuum~$X$ and a lattice-base~$\famB$ for its closed sets be given.
As before we start with the formulas that ensure that $wL$~will be a 
hereditarily indecomposable continuum.
To these formulas we add the diagram of~$\famB$; this consists of~$\famB$
itself, as a set of constants, and the `multiplication tables' for~$\meet$
and~$\join$, i.e., $A\meet B=C$ whenever $A\cap B=C$ and $A\join B=C$
whenever $A\cup B=C$.

Now, if $L$~is to satisfy the diagram of~$\famB$ it must contain 
elements~$x_A$ for every $A\in\famB$ so that $x_{A\meet B}=x_A\meet x_B$
and $x_{A\join B}=x_A\join x_B$ hold whenever appropriate; but this simply
says that there is an embedding of $\famB$ into~$L$.

We are left with the task of showing that the set~$\famT$ of formulas that 
express distributivity, disjunctiveness, normality as well as
formulas~\ref{L.is.connected} and~\ref{L.is.HI} (and~\ref{dimL<=1}), 
together with the diagram of~$\famB$ is consistent.
Let $T$ be a finite subset of~$\famT$ and, if necessary, add the first six 
formulas to it.
Let $\famB'$ be a countable, normal and disjunctive sublattice of~$\famB$
that contains the finitely many constants that occur in~$T$.
The Wallman space of~$\famB'$, call it~$Y$, is a metric continuum
and therefore the continuous image of a hereditarily indecomposable
(one-dimensional) continuum~$K$.
The lattice of closed sets of~$K$ satisfies all the formulas from~$T$:
interpret $A$ by its preimage in~$K$.

It follows that $\famT$~is consistent and that it therefore has a model~$L$
of the same cardinality as~$\famT$, which is the same as the cardinality
of~$\famB$.
The lattice~$L$ satisfies all formulas from~$\famT$; its Wallman space
is a (one-dimensional) hereditarily indecomposable continuum that maps
onto~$X$.
If $\famB$ is chosen to be of minimal size then $wL$ is of the same weight
as~$X$.

\subsubsection{Making a weakly confluent map}

We now improve the foregoing construction so as to make the continuous
surjection weakly confluent.

The following theorem --- which is a souped-up version of the Marde\v{s}i\'c
factorization theorem --- implies that it suffices to get some hereditarily 
indecomposable continuum~$Y$ that admits a weakly confluent map~$f$ 
onto our continuum~$X$.

\begin{theorem}
Let $f:Y\to X$ be a continuous surjection between compact Hausdorff spaces. 
Then $f$ can be factored as~$h\circ g$, where 
$Y\buildrel g\over\to Z\buildrel h\over\to X$ and $Z$ has the same weight 
as~$X$ and shares many properties with~$Y$.
\end{theorem}

\begin{proof}
Let $\famB$ be a lattice-base for the closed sets of~$X$ (of minimal size) 
and identify it with its copy $\{f\preim[B]:B\in\famB\}$ in~$\Hyp{Y}$.
By the L\"owenheim-Skolem theorem 
\cite[Corollary 3.1.5]{Hodges93} there is a lattice~$\famD$, of the same
cardinality as~$\famB$, such that $\famB\subseteq\famD\subseteq\Hyp{Y}$ 
and $\famD$~is an elementary substructure of~$\Hyp{Y}$.
The space~$Z=w\famD$ is as required.
\end{proof}

Some comments on this theorem and its proof are in order, because they do
not seem to say very much.
However, at this point we can see the power of the notion of an
elementary substructure.
From the knowledge that the smaller lattice~$\famD$ contains
solutions for every equation with parameters from~$\famD$ that is
solvable in~$\Hyp{Y}$ we can deduce a lot about~$Z$.

For example, if $Y$ hereditarily indecomposable then so is~$Z$.
For if $(C,D,F,G)$~is a pliand foursome in~$\famA$ then the equation
\begin{multline}
(C\meet(z_2\join z_3)=\0)\land (D\meet(z_1\join z_2)=\0)
   \land (z_1\meet z_3=\0){}\\
{}\land(z_1\meet z_2\meet G=\0)\land(z_2\meet z_3\meet F=\0)\land
 (z_1\join z_2\join z_3=\1)    
\end{multline}
has a solution in~$\Hyp{Y}$, hence in~$\famD$.

A similar argument establishes $\dim Z=\dim Y$: the Theorem on 
Partitions
(\cite[Theorem~7.2.15]{Engelking89}) yields systems of equations
that characterize covering dimension.
For example, if $\dim Y\le1$ then $\dim Z\le 1$ because
if $A,B,C,D\in\famA$ satisfy $A\cap B=C\cap D=\emptyset$ then
$\Hyp{Y}$, and hence~$\famA$, contains a solution to 
\begin{multline}
(A\meet u_0=\0)\land (B\meet v_0=\0)\land
                (C\meet u_1=\0)\land (D\meet v_1=\0)\land{}\\
        {}\land (u_0\join v_0=\1)\land (u_1\join v_1=\1)\land
                (u_0\meet v_0\meet u_1\meet v_1=\0).
\end{multline}
The negation of formula~\ref{dimL<=1} yields a (parameterless)
equation that has a solution in~$\Hyp{Y}$ iff $\dim Y>1$.
We invite the reader to explore how the solution that must exist
in~$\famA$ witnesses that $\dim Z>1$.

We leave to the reader the verification that if $f$ is weakly confluent
then so is the map~$h$ in the factorization.

Now let $X$ be a continuum.
Our aim is of course to find a lattice~$L$ that contains the diagram 
of~$\HypX$ --- to get our continuous surjection~$f$ --- 
and for every~$C\in\CsX$ a continuum~$C'$ in~$wL$ such that $f[C']=C$; here $\CsX$~denotes the family of subcontinua of~$X$.

As before we add the diagram of~$\HypX$ to the formulas that guarantee
that $wL$ will be a hereditarily indecomposable continuum.
In addition we take a set of constants $\{C':C\in\CsX\}$ and stipulate
that $C'$ will be a continuum that gets mapped onto~$C$.

To make sure that every $C'$ is connected we put $\conn(C')$ into
our set of formulas, for every~$C$.
Next, $f[C']\subseteq C$ translates, via the embedding into~$L$, into
$C'\le C$ (or better $C'=C'\meet C$).
Now, if it happens that $f[C']\subsetneq C$ then there is a closed set~$D$
in~$X$ (in fact it is $f[C']$ but that is immaterial) such that
$C'\le D$ and $C\notle D$.
In order to avoid this we also add, for every~$C\in\CsX$ and 
every~$D\in\HypX$, the formula
$$
(C'\le D) \implies (C\le D)
$$
to our set of formulas.

Again, the theorem in the metric case implies that this set of formulas is
consistent --- given a finite subset~$T$ of it make a metric 
continuum~$X_T$ as before, by expanding $\{B\in\HypX:B$~occurs in~$T\}$
to a countable normal sublattice~$\famB$ of~$\HypX$; 
then find a metric continuum~$Y_T$ of the desired type that admits a weakly
confluent map~$f$ onto~$X_T$; finally choose for every 
$C\in\CsX$ that occurs in~$T$ a continuum in~$Y_T$ that maps onto~$C$ and 
assign it to~$C'$; this then makes the family of closed sets of~$Y_T$ a model
of~$T$.

As before we obtain a lattice~$L$ whose Wallman space is one-di\-men\-sional
and hereditarily indecomposable, and which, in addition, admits a weakly confluent
map onto~$X$.

\section{$\OCA$ and some of its uses}
\label{sec.OCA}

The \emph{Open Coloring Axiom} ($\OCA$) was formulated by
To\-dor\-\v{c}e\-vi\'c in~\cite{Todorcevic89}.
It reads as follows:
if $X$ is separable and metrizable and if $[X]^2=K_0\cup K_1$,
where $K_0$~is open in the product topology of~$[X]^2$, then
\emph{either} $X$~has an uncountable $K_0$-homogeneous subset~$Y$
\emph{or} $X$~is the union of countably many $K_1$-homogeneous subsets.

One can deduce $\OCA$ from the \emph{Proper Forcing Axiom} ($\PFA$)
or prove it consistent in an $\omega_2$-length countable support
proper iterated forcing construction, using~$\diamond$ on~$\omega_2$
to predict all possible subsets of the Hilbert cube and all possible
open colourings of these.

The axiom $\OCA$ has a strong influence on the structure of maps between
concrete objects like $\powN$, $\powNfin$ and the measure algebra~$\Meas$.
In fact it imposes such strict conditions that $\OCA$ implies the 
nonembeddability of~$\Meas$ and other algebras into~$\powNfin$.

\subsection{A simple space}\label{subsec.Dstar}

Let $\D=\omega\times(\omega+1)$; Parovi\v{c}enko's theorem implies
that $\Dstar$~is a continuous image of~$\Nstar$ if $\CH$~is assumed.
We shall see that such a continuous surjection has no simple description.
Later on we shall indicate how $\OCA$~dictates that if there is a continuous
surjection of~$\Nstar$ onto~$\Dstar$ at all then there must also be one
with a simple description, thus showing that $\OCA$~implies
$\Dstar$~is \emph{not} a continuous image of~$\Nstar$.

Most of the proof will be algebraic, i.e., instead of working with
continuous maps from~$\Nstar$ onto~$\Dstar$ we work with embeddings
of the algebra of clopen sets of~$\Dstar$ into~$\powNfin$.
However, both algebras are quotient algebras so we will consider
liftings of these embeddings, i.e., we will work with maps
from~$\Co\D$ to~$\powN$ that represent them.

First of all we give a description of the Boolean algebra of clopen
subsets of~$\D$ that is easy to work with.
We work in $\omega\times\omega$ and denote the $n$-th column
$\{n\}\times\omega$ by~$C_n$.
The family
$$
\Ba{B}=\bigl\{X\subseteq\omega\times\omega:
             (\forall n\in\omega)
             (C_n\almsubseteq X \lor C_n\cap X\almeq\emptyset)\bigr\}
$$
is the Boolean algebra of clopen subsets of~$\D$.
We also consider the subfamily
$$
\Ba{B}^-=\bigl\{X\in\Ba{B}:
              (\forall n\in\omega)
              (C_n\cap X\almeq\emptyset)\bigr\}
$$
of~$\Ba{B}$.

Now assume $S:\Nstar\to\Dstar$ is a continuous surjection and take
a map~$\Sigma:\Ba{B}\to\powN$ that represents~$S$,
i.e., for all $X\in\Ba{B}$ we have $\Sigma(X)^*=S\preim[X^*]$.
Note that if $X$ is compact in~$\D$ then $\Sigma(X)$ is finite.

The main result of this section is that $\Sigma$ cannot be 
\emph{simple}, where simple maps are defined as follows.

\begin{definition}
We call a map $F:\Ba{B}^-\to\pow(\omega)$ \emph{simple}
if there is a map $f$ from $\omega\times\omega$ to $\finsets$
such that
$F(X)\almeq f[X]$ for all~$X$, where $f[X]$ denotes the set
$\bigcup_{x\in X}f(x)$.
\end{definition}

\begin{theorem}\label{no.simple.maps}
The map $\Sigma\restr\Ba{B}^-$ is not simple.
\end{theorem}

\begin{proof}
We assume that there is a
map~$\sigma:\omega\times\omega\to\finsets$
such that $\sigma[X]\almeq\Sigma(X)$ for all~$X$;
this implies that $\sigma[X]^*=S\preim[X^*]$ for all~$X$, so
the map $X\mapsto\sigma[X]$ also represents~$S$.
We may therefore as well assume that $\Sigma(X)=\sigma[X]$ for all~$X$.

\begin{claim}
We can assume that the values  $\sigma(x)$ are pairwise disjoint.
\proof
Let $\langle f_\alpha:\alpha<\bee\rangle$ be a sequence
in~$\functions$ that is strictly increasing and unbounded with respect
to~$\alml$; also each $f_\alpha$ is assumed to be strictly increasing.

For each $\alpha$ let $L_\alpha=\bigl\{(n,m):m\le f_\alpha(n)\bigr\}$ and
$A_\alpha=\sigma[L_\alpha]$.
Next let
$$
B_\alpha=\bigr\{i\in A_\alpha:
          (\exists x,y\in L_\alpha)
     \bigl((x\neq y)\land (i\in\sigma(x)\cap\sigma(y))\bigr)\bigr\}.
$$
Now if $B_\alpha$ were infinite then we could find different~$i_n$
in~$B_\alpha$ and different $x_n$ and $y_n$ in~$L_\alpha$ such that
$i_n\in\sigma(x_n)\cap\sigma(y_n)$.
But then $X=\{x_n:n\in\omega\}$ and $Y=\{y_n:n\in\omega\}$ would be
disjoint yet $\sigma[X]\cap\sigma[Y]$ would be infinite.

We conclude that each $B_\alpha$ is finite and because $\bee$~is regular
we can assume that all~$B_\alpha$ are equal to the same set~$B$.
Fix~$n$ such that $[n,\omega)\times\omega\subseteq\bigcup_\alpha L_\alpha$
and note that on $[n,\omega)\times\omega$ we have
$\sigma(x)\cap\sigma(y)\subseteq B$ whenever $x\neq y$.
Replace $\sigma(x)$ by $\sigma(x)\setminus B$ and $\omega\times\omega$
by $[n,\infty)\times\omega$.
\end{claim}

In a similar fashion we can prove the following claim.

\begin{claim}
We can assume that the values $\sigma(x)$ are all nonempty.
\proof
There are only finitely many $n$ for which there is an~$m$ such that
$\sigma(n,m)=\emptyset$.
Otherwise we could find a noncompact~$X\in\Ba{B}^-$ for
which~$\Sigma(X)=\emptyset$.
Drop these finitely many columns from $\omega\times\omega$.
\end{claim}

For each $n$ let $D_n=\sigma[C_n]$ and work inside~$D=\bigcup_nD_n$.
Also define, for $f\in\functions$, the sets
$L_f=\bigl\{(n,m): m\le f(n)\bigr\}$ and $M_f=\sigma[L_f]$.

Now observe the following:
for each $f$ and $n$ the intersection $M_f\cap D_n$ is finite and
if $X\subseteq D$ is such that $X\cap D_n\almeq\emptyset$ for all~$n$ then
$X\subseteq M_f$ for some~$f$.

In $\Dstar$ we consider the top line
$T=\bigl(\omega\times\{\omega\}\bigr)^*$
and its complement~$O$.
First we note that $O=\bigcup_fL_f^*$ and so
$$
S\preim[O]=\bigcup_fS\preim[L_f^*]=\bigcup_f\sigma[L_f]^*=\bigcup_fM_f^*.
$$
This means that $S[D_n^*]\subseteq T$ for all~$n$, because $D_n^*$~is
disjoint from $\bigcup_fM_f^*$.
Also, the boundary of the cozero set $\bigcup_nD_n^*$ is the boundary
of~$\bigcup_fM_f^*$; by continuity this boundary is mapped onto the
boundary of~$O$, which is~$T$.

This argument works for every infinite subset~$A$ of~$\omega$:
the boundary of $\bigcup_{n\in A}D_n^*$ is mapped exactly onto
the set~$T_A=\bigl(A\times\{\omega\}\bigr)^*$ and so
$T_A$~is contained in the closure of $S[\bigcup_{n\in A}D_n^*]$
and $S[D_n^*]\subseteq T_A$ for all but finitely many~$n\in A$.

From the fact that nonempty $G_\delta$-sets in~$\Nstar$ have nonempty
interior one readily deduces that no countable family of nowhere dense
subsets of~$\Nstar$ has a dense union.
We conclude that there is an~$n_0$ such that $\Int_TS[D_{n_0}^*]$ is
nonempty.
Choose an infinite subset $A_0$ of~$(n_0,\omega)$ such that
$T_{A_0}\subseteq S[D_{n_0}^*]$.

Continue this process: once $n_i$ and $A_i$ are found one finds
$n_{i+1}\in A_i$ such that $S[D_{n_{i+1}}^*]$ has nonempty interior and
is contained in~$T_{A_i}$, next choose an infinite subset~$A_{i+1}$
of~$A_i\cap(n_{i+1},\omega)$ such that
$T_{A_{i+1}}\subseteq S[D_{n_{i+1}}^*]$.

Finally then let $A=\{n_{2i}:i\in\omega\}$ and $B=\{n_{2i+1}:i\in\omega\}$.
Note that $T_A\subseteq \bigcap_{n\in B}S[D_n^*]$
but also that $S[D_n^*]\cap T_A=\emptyset$ for all but finitely
many~$n\in B$.
This contradiction completes the proof of
Theorem~\ref{no.simple.maps}.
\end{proof}

\subsection{The measure algebra}\label{subsec.meas}

Parovi\v{c}enko's theorem also implies that the measure algebra~$\Meas$
can be embedded into~$\powNfin$ (under $\CH$); as in the previous section
we shall see that such an embedding admits no easy description.
Again, $\OCA$~dictates that any embedding of~$\Meas$ into~$\powNfin$ induces
an embedding with an easy description, from
which we deduce that $\OCA$~prohibits embeddability of~$\Meas$ into~$\powNfin$.

\subsubsection{The Measure Algebra}

The standard representation of the Measure Algebra is as the quotient of
the $\sigma$-algebra of Borel sets of the unit interval by the ideal
of sets of Lebesgue measure zero.
For ease of notation we choose a different underlying set, namely
$\C=\omega\times\Cantor$, where $\Cantor$is the Cantor set.
We consider the Cantor set endowed with the natural coin-tossing
measure~$\mu$, determined by specifying $\mu\bigl([s]\bigr)=2^{-\card{s}}$.
Here $s$ denotes a finite partial function from~$\omega$ to~$2$
and $[s]=\{x\in\Cantor:s\subset x\}$.
We extend $\mu$ on the Borel sets of~$\C$ by setting
$\mu\bigl(\{n\}\times[s]\bigr)=2^{-\card{s}}$ for all~$n$ and~$s$.

The measure algebra is isomorphic to the quotient algebra
$\Meas=\Bor(\C)/\Null$, where $\Null=\{N\subseteq\C:\mu(N)=0\}$;
henceforth we shall work with $\Meas$.

\subsubsection{Liftings of embeddings}

Assume $\phi:\Meas\to\powNfin$ is an embedding of Boolean algebras and take a
\emph{lifting} $\Phi:\Meas\to\powN$ of~$\phi$; this is a map that chooses a
representative~$\Phi(a)$ of~$\phi(a)$ for every~$a$ in~$\Meas$.

We shall be working mostly with the restrictions of~$\phi$ and~$\Phi$ to the
family of (equivalence classes of) open subsets of~$\C$ and in particular with
their restrictions to the canonical base for~$\C$, which is
$$
\famB=\bigl\{\{n\}\times[s]:n\in\omega, s\in\bintree\bigr\}.
$$
To keep our formulas manageable we shall identify $\famB$ with the
set $\omega\times\bintree$.
We shall also be using layers/strata of~$\famB$ along functions from~$\omega$
to~$\omega$: for $f\in\functions$ we put
$\famB_f=\bigl\{\orpr{n}{s}:n\in\omega, s\in\sbintree\bigr\}$.

For a subset $O$ of~$\famB_f$ we abbreviate
$\phi\bigl(\bigcup\bigl\{\{n\}\times[s]:\orpr{n}{s}\in O\bigr\}\bigr)$
by~$\phi(O)$ and define $\Phi(O)$ similarly.
Observe that $O\mapsto\phi(O)$ defines an embedding of~$\pow(\famB_f)$
into~$\powNfin$.
As an extra piece of notation we use $\Phi[O]$ (square brackets) to denote the
union $\bigcup\bigl\{\Phi(n,s):\orpr{n}{s}\in O\bigr\}$,
where $\Phi(n,s)$ abbreviates $\Phi\bigl(\bigl\{\orpr{n}{s}\bigr\}\bigl)$.

For later use we explicitly record the following easy lemma.

\begin{lemma}\label{lemma.Phi.fin.add}
If $f\in\functions$ and if\/ $O$ is a finite subset of\/~$B_f$ then
$\Phi(O)\almeq\Phi[O]$.
\end{lemma}

\begin{proof}
Both sets represent $\phi(O)$.
\end{proof}

Let us call a lifting \emph{complete} if it satisfies
Lemma~\ref{lemma.Phi.fin.add} for \emph{every}~$f\in\functions$ and
\emph{every} subset~$O$ of~$B_f$.

We can always make a lifting~$\Phi$ \emph{exact}, by which we mean that
the sets~$\Phi(n,\emptyset)$ form a partition of~$\N$ and that
every~$\Phi(s,n)$ is the disjoint union of~$\Phi(n,\szer)$
and~$\Phi(n,\sone)$;
indeed, we need only change each of the countably many sets~$\Phi(n,s)$ by
adding or deleting finitely many points to achieve this.

Now we can properly formulate what `easy description' means and how $\OCA$
insists on there being an easily described embedding.
\begin{enumerate}
\item For every exact lifting~$\Phi$ of an embedding~$\phi$ there are
      an~$f\in\functions$ and an infinite subset~$O$ of~$B_f$ such that
      $\Phi(O)\almneq\Phi[O]$, i.e., no exact lifting is complete
      --- see Proposition~\ref{prop.exact.is.not.complete}.
\item $\OCA$~implies that every embedding~$\phi$ gives rise to an
      embedding~$\psi$ with a lifting~$\Psi$ that is both exact and
      complete (see \cite{DowHart2000a}).
\end{enumerate}

\subsubsection{No exact lifting is complete}
\label{sec.exact.is.not.complete}

Assume $\phi:\Meas\to\powNfin$ is an embedding and consider an exact lifting~$\Phi$
of~$\phi$.
The following proposition shows that $\Phi$~is not complete.

\begin{proposition}\label{prop.exact.is.not.complete}
There is a sequence $\langle t_n:n\in\omega\rangle$ in~$\bintree$
such that for the open set $O=\bigcup_{n\in\omega}\{n\}\times[t_n]$ we have
$\Phi(O)\almneq\Phi[O]$.
\end{proposition}

\begin{proof}
Take, for each~$n$, the monotone enumeration~$\{k(n,i):i\in\omega\}$
of~$\Phi(n,\emptyset)$ and apply exactness to find $t(n,i)\in2^{i+2}$
such that $k(n,i)\in\Phi\bigl(n,t(n,i)\bigr)$.
Use these~$t(n,i)$ to define open sets
$U_n=\bigcup_{i\in\omega}\{n\}\times\bigl[t(n,i)\bigr]$;
observe that $\mu(U_n)\le\sum_{i\in\omega}2^{-i-2}=\frac12$.
It follows that $\Phi\bigl(\{n\}\times U_n^c\bigr)$ is infinite.

We let $F$ be the closed set $\bigcup_{n\in\omega}\{n\}\times U_n^c$;
its image~$\Phi(F)$ meets every~$\Phi(n,\emptyset)$ in an infinite set.
For every~$n$ let $i_n$ be the first index with $k(n,i_n)\in\Phi(F)$
and consider the open set
$O=\bigcup_{n\in\omega}\{n\}\times\bigl[t(n,i_n)\bigr]$
and the infinite set~$I=\bigl\{k(n,i_n):n\in\omega\bigr\}$.

Observe the following
\begin{enumerate}
\item $\Phi(O)\cap\Phi(F)\almeq\emptyset$, because $O\cap F=\emptyset$;
\item $I\subseteq\Phi(F)$, by our choice of the $i_n$; and
\item $I\subseteq\Phi[O]$, by the choice of the $t(n,i_n)$.
\end{enumerate}
It follows that $\bigl<t(n,i_n):n\in\omega\bigr>$ is as required.
\end{proof}

\subsection{How $\OCA$ induces simple structure}

In the previous two subsections we had two maps,
$\Sigma:\Ba{B}^-\to\powN$ and $\Phi:\famB\to\powN$.
Both induced embeddings of their domains into~$\powNfin$.
What $\OCA$ does is guarantee the existence of an infinite subset~$A$
of~$\omega$ such that $\Sigma$ is simple on 
$\{B\in\Ba{B}^-:B\subseteq A\times\omega\}$ and
such that the embedding induced by~$\Phi$ has a lifting that is exact and
complete on $\{\orpr ns\in\famB:n\in A\}$.
We indicate how to do this for~$\Sigma$ and refer the interested reader to
\cite{DowHart2000a} for details on how to deal with~$\Phi$.

\subsubsection{Working locally}

For $f\in\functions$ and put $L_f=\{\orpr mn:n\le f(m)\}$ and observe
that for every~$B\in\Ba{B}^-$ there is an~$f$ such that $B\subseteq L_f$;
this means that $\Ba{B}^-=\bigcup_f\pow(L_f)$.

Our first step, for $\Sigma$, will be to show that it is simple on 
$L_{f,A}=\orpr mn\in L_f:m\in A\}$ for many subsets of~$\omega$ (for all~$f$).
Similarly, for $\Phi$, we show that there is an exact and complete lifting
on~$\famB_{f,A}=\{\orpr ns\in\famB_f:n\in A\}$ for many subsets of~$\omega$
(for all~$f$).
The proof will be finished by finding one~$A$ that works for all~$f$ 
simultaneously.
We follow the strategy laid out in Veli\v{c}kovi\'c' 
papers \cite{Velickovic86} and \cite{Velickovic93}.

Fix a bijection $c:\omega\to\bintree$ and use it to transfer the set of
branches to an almost disjoint family~$\famA$ on~$\omega$ and
fix an $\aleph_1$-sized subfamily $\{A_\alpha:\alpha<\omega_1\}$ of~$\famA$.
Using~$\OCA$ we shall show that all but countably many~$A_\alpha$ are as
required, i.e, $\Phi$~is simple on $L_{f,A_\alpha}$ for all but countably 
many~$\alpha$.
Let us write $L_\alpha=L_{f,A_\alpha}$.

To apply $\OCA$ we need a separable metric space; we take
$$
X=\bigl\{\orpr{a}{b}:
  (\exists\alpha<\omega_1)(b\subseteq a\subseteq L_\alpha)\bigr\},
$$
topologized by identifying $\orpr{a}{b}$
with~$\bigl\langle a,b,\Sigma(a),\Sigma(b)\bigr\rangle$ ---
that is, $X$~is identified with a subset of~$\pow(\omega)^4$.
We define a partition $[X]^2=K_0\cup K_1$ by:
$\bigl\{\orpr{a}{b},\orpr{c}{d}\bigr\}\in K_0$ \iff/
1)~$a$ and~$c$ are in different~$L_\alpha$'s;
\ 2)~$a\cap d=c\cap b$, and
\ 3)~$\Sigma(a)\cap\Sigma(d)\neq \Sigma(c)\cap\Sigma(b)$.

Because of the special choice of the almost disjoint family~$\famA$ 
the set~$K_0$ is open: condition~1) can now be met using only finitely many
restrictions and then condition~2) needs finitely many restrictions also;
condition~3) needs just one restriction.

The next step is to show that there is no uncountable $K_0$-homogeneous set.
Suppose $Y$ were uncountable and $K_0$-homogeneous.
Then we can form the set $x=\bigcup\{b:(\exists a)(\orpr ab\in Y)\}$.
Condition~2) implies that $x\cap a=b$ whenever $\orpr ab\in Y$ and
this means that $\Sigma(x)\cap\Sigma(a)\almeq\Sigma(b)$ for all these
pairs.
So now we can fix $n\in\omega$ and subsets $p$ and~$q$ of~$n$ such that,
for uncountably many $\orpr ab\in Y$ we have
$\bigl(\Sigma(x)\cap\Sigma(a)\bigr)\symmdif\Sigma(b)\subseteq n$,
$\Sigma(a)\cap n=p$ and $\Sigma(b)\cap n=q$.
But then condition~3) would be violated for these pairs.

We conclude that $X=\bigcup_nX_n$, where each $X_n$ is $K_1$-homogeneous.
Choose, for each~$n$, a countable dense set~$D_n$ in~$X_n$ --- with respect
to the given topology.
Let $\alpha_f$ be the first ordinal such that if $\orpr ab\in\bigcup_nD_n$
and $a\subseteq L_\alpha$ then $\alpha<\alpha_f$.
For $\alpha\ge\alpha_f$ and $n\in\omega$ define $F_n:\pow(L_\alpha)\to\powN$
by
$$
F_n(b)=\bigcup\{\Sigma(L_\alpha)\cap\Sigma(d):
                (\exists c)(\orpr cd\in D_n \land c\cap b=L_\alpha\cap d)\}.
$$
Each of the maps $F_n$ is Borel and 
$F_n(b)=\Sigma(b)$ whenever $\orpr {L_\alpha}{b}\in X_n$.
Thus $\Sigma$ has been tamed substantially: it has been covered by countably
many Borel maps.
In~\cite{DowHart99} one can find how to modify Veli\v{c}kovi\'c' arguments
from~\cite{Velickovic93} to show that this implies that $\Sigma$~is indeed
simple.

\subsubsection{Going global}

We now have for each~$f$ an ordinal~$\alpha_f$ such that $\Sigma$~is simple
on~$L_{f,A_\alpha}$ whenever $\alpha\ge\alpha_f$.
It should be clear that in case $f\almle g$ and $\Sigma$~is simple 
on~$L_{g,A_\alpha}$ it is also simple on~$L_{f,A_\alpha}$ because the latter
set is almost a subset of the former.
It follows that $f\mapsto\alpha_f$ is monotone from $\functions$ to~$\omega_1$.

Now, $\OCA$ implies that $\bee=\aleph_2$, see~\cite[Theorem~3.16]{Bekkali91}.
But this then implies that there is an ordinal~$\alpha_\infty$ such that
$\alpha_f\le\alpha_\infty$ for all~$f$.
We find that, for very $\alpha\ge\alpha_\infty$, the map~$\Sigma$ is simple
on~$L_{f,A_\alpha}$ for all~$f$.

For definiteness let $A=A_{\alpha_\infty}$ and fix for each~$f$ a map
$\sigma_f:L_{f,A}\to\finsets$ that induces~$\Sigma$.
It should be clear that the $\sigma_f$ cannot differ too much, i.e,
on $L_{f,A}\cap L_{g,A}$ the maps $\sigma_f$ and $\sigma_g$ will differ
in only finitely many point --- the family $\{\sigma_f:f\in\functions\}$
is said to be coherent.
Theorem~3.13 from~\cite{Bekkali91} now applies: one can find one map
$\sigma:A\times\omega\to\finsets$ such that 
$\sigma\restr L_{f,A}\almeq\sigma_f$ for all~$f$.
This $\sigma$~is the simplifying map that we were looking for.

\let\Pow\pow
\providecommand{\cprime}{$'\!$}
\providecommand{\cyr}{}
\bibliographystyle{amsplain}
\bibliography{names,longnames,research}

\end{document}